\documentclass{amsart}
\usepackage{amsmath,amsthm,amssymb}
\usepackage[unicode]{hyperref}
\hypersetup{
	colorlinks=true,
	linkcolor=black,
	citecolor=black,
	urlcolor=blue
}
\usepackage{enumitem}

\newtheorem{thm}{Theorem}[section]

\newtheorem{conj}[thm]{Conjecture}
\newtheorem{cor}[thm]{Corollary}
\newtheorem{claim}[thm]{Claim}
\theoremstyle{definition}
\newtheorem{defn}[thm]{Definition}

\numberwithin{equation}{section}

\newcommand{\infsub}[1]{{[#1]^\omega}}
\newcommand{\idl}[1]{{\mathcal #1}}
\newcommand{\zfaxiom}{{\sf ZF}}

\newenvironment{claimproof}[1]{%
	\renewcommand{\qedsymbol}{{\tiny Claim #1} \(\square\)}%
	\begin{proof}[Proof of Claim #1]%
	}{%
	\end{proof}%
}

\title{A proof of Mathias's conjecture}

\author{Jialiang He}
\address{College of Mathematics, Sichuan University, Chengdu, Sichuan, 610064 China}
\email{jialianghe@scu.edu.cn}

\author[Jintao Luo]{Jintao Luo$ ^\ast $}
\address{College of Mathematics, Sichuan University, Chengdu, Sichuan, 610064 China}
\email{jintaoluo@foxmail.com}

\author{Shuguo Zhang}
\address{College of Mathematics, Sichuan University, Chengdu, Sichuan, 610064 China}
\email{zhangsg@scu.edu.cn}



\thanks{$^{\ast}$Corresponding author}
\begin{document}
	
	\begin{abstract}
		We prove in $\mathsf{ZF}$ that the generalized Ramsey theorem $\omega\longrightarrow(\omega)^\omega_2$ implies the nonexistence of infinite maximal almost disjoint families, solving a conjecture of A.~R.~D.~Mathias from 1977.
	\end{abstract}
	
	\maketitle
	
	\section{Introduction}
	
	Let $ \omega $ be the set of non-negative integers. Ramsey's theorem $\omega\longrightarrow(\omega)^n_2$ states that for each positive integer $n$, every two-coloring of the $n$-element subsets of $\omega$ admits an infinite homogeneous set. What happens when $n$ is replaced by $\omega$?

Denote the collection of infinite subsets of $ \omega $ by $ [\omega]^\omega $. Identifying each subset $X\subseteq[\omega]^\omega$ with its corresponding two-coloring, we say that $X$ is \emph{Ramsey} if this coloring admits an infinite homogeneous set. The generalized Ramsey theorem $\omega\longrightarrow(\omega)^\omega_2$ is then equivalent to the statement that every subset of $[\omega]^\omega$ is Ramsey. Erd{\H{o}}s and Rado~\cite{erdosrado1952classifications} showed that there is a non-Ramsey set under the \emph{Axiom of Choice}. As their argument is nonconstructive, it is natural to ask which definable sets have the Ramsey property. Galvin and Prikry~\cite{galvinprikry1973borel} proved that all Borel sets have this property, and Silver~\cite{silver1970ramsey} extended this to analytic sets. In his doctoral thesis, Mathias~\cite{mathias1977happy} established the relative consistency of the Ramsey property for all sets in Solovay's model.

Based on the behavior of non-Ramsey sets under various conditions listed above, we can ask similar questions for other combinatorial objects constructed using choice. How simply can we define such objects? Do they exist in Solovay’s model?

Regarding these questions, Mathias also studied several such objects in his 1977 paper~\cite{mathias1977happy}, including infinite maximal almost disjoint (MAD) families. A family $\mathcal A\subseteq[\omega]^\omega$ is \emph{almost disjoint} if any two distinct members have finite intersection. Such a family is \emph{maximal} if every infinite subset of $\omega$ has infinite intersection with some member of $\mathcal A$. Mathias used Ramsey-theoretic methods to prove there are no analytic infinite MAD families~\cite[Corollary~4.7]{mathias1977happy}. Mathias also showed that, under a further large cardinal assumption, there are no infinite MAD families in Solovay's model. This left open the following two conjectures.
	
	\begin{conj}[\cite{mathias1977happy}]{\label{conjecture:con_ramsey}}
		\phantom{a}
		\begin{enumerate}[label=(\arabic*), font=\upshape]
			\item Are there infinite MAD families in Solovay's model?
			\item \cite[p.~87]{mathias1977happy} Is it a theorem of $\zfaxiom$ and $ \mathsf{DC}_{\mathbb R} $ that there are no infinite MAD families assuming all sets are Ramsey?
		\end{enumerate}
	\end{conj}
	
	Here $\mathsf{DC}_{\mathbb R}$ denotes the \emph{Axiom of Dependent Choice for Reals}.
	
	The first conjecture was answered negatively by T\"ornquist in~\cite{tornquist2018definability}. Horowitz and Shelah~\cite{horowitz2019non} settled the question of consistency strength by showing that the nonexistence of infinite MAD families is equiconsistent with $\mathsf{ZFC}$. Neeman and Norwood~\cite{neeman2018happy} strengthened T\"ornquist's result and showed that there are no infinite MAD families under a variant of the \emph{Axiom of Determinacy} $\mathsf{AD}^+$. Bakke Haga, Schrittesser and T\"ornquist~\cite{bakke2022maximal} further extended these results to higher dimensional MAD families. Recently, the results on determinacy and MAD families were generalized to uncountable cardinals by Chan, Jackson and Trang~\cite{chan2024almost}.
	
	For the second conjecture, Schrittesser and T\"ornquist~\cite{schrittesser2019ramsey} proved that if all sets have the Ramsey property, then there are no infinite MAD families assuming \emph{Dependent Choice} $ \mathsf{DC} $ and \emph{Ramsey Uniformization} $ \mathsf{R}$-$\mathsf{Unif} $. They subsequently extended this result to MAD families with respect to all transfinitely iterated Fr\'echet ideals~\cite{schrittesser2020ramsey}, under similar choice principles.
	
	In this paper, we prove a stronger form of Conjecture~\ref{conjecture:con_ramsey}(2) without any choice assumptions, and recover several earlier results in the literature. To state our main theorem, recall a \emph{boldface pointclass} is a collection $\Gamma$ of subsets of Polish spaces which is closed under continuous preimages and containing all clopen sets \cite[p.~118]{kechris2006axiom}. For our purposes, we say a pointclass is \emph{good} if it is moreover closed under continuous images and bounded universal number quantification. Each $\boldsymbol{\Sigma}^1_n$ for $ n\geq 1 $ is a good pointclass (cf.~\cite[Theorem~1E.2]{moschovakis2025descriptive}).

	\begin{thm}[$\zfaxiom$]\label{thm:mad-ramsey}
		Let $\Gamma$ be a good pointclass. If every subset of $ \infsub{\omega} $ in $\Gamma$ has the Ramsey property, then there are no infinite MAD families in $\Gamma$.
	\end{thm}
	
	Conjecture~\ref{conjecture:con_ramsey}(2) follows in $\zfaxiom$ alone by taking $\Gamma$ to be the pointclass of all sets. The theorem also directly gives short proofs of several results in the literature, showing the nonexistence of infinite MAD families within the pointclass of analytic sets, Solovay's model and under $ \mathsf{AD}^+$~\cite{mathias1977happy,tornquist2018definability,neeman2018happy}.
	
	The key idea of our proof is a \emph{hiding} process. In the absence of choice principles, maximality only provides a multivalued total relation. However, through a suitable encoding operator, one can shrink into the infinite intersection and densely find an encoded real contained in a \emph{single} MAD member. To further provide a non-Ramsey set, we also need a mechanism to create local oscillation, and then densely escape from hiding.
	
	In the Dedekind-infinite case, the tilde operator of Schrittesser and T\"ornquist ~\cite{schrittesser2019ramsey} provides a natural escaping mechanism; the challenge is to make hiding compatible. Inspired by an earlier construction of Mathias~\cite{mathias1975remark}, we define a new encoding operator such that both hiding and splitting of encoded reals occur densely.
	
	In the Dedekind-finite case, if we only have the trivial hiding process, oscillation is impossible, since hiding is preserved under shrinking. However, Dedekind-finiteness gives a oscillatable hiding process. Each real can be encoded via iterated hidings to determine distinct MAD members, which terminates in finite steps; while maximality densely extends it to arbitrarily long lengths. This provides the desired hiding and oscillation mechanism.
	
	\section{Proof of the main theorem}
	
	For any set $X$, let $[X]^\omega$ be the set of infinite subsets of $X$. We equip $[\omega]^\omega$ with the subspace topology inherited from $2^\omega$ via characteristic functions. For any $ n\in\omega $, let $ X^n $ be the set of sequences in $ X $ with length $ n $, and let $ X^{<\omega} $ be the set of finite sequences in $ X $. If $ X $ is a Polish space, we equip $ X^{<\omega} $ with the topological sum. For $x\in[\omega]^\omega$, we write $x=\{x(0),x(1),\ldots\}$, where $x(0)<x(1)<\cdots$. For $x,y\subseteq\omega$, write $x\subseteq^*y$ if $x\setminus y$ is finite. Given a function $ f:X\to Y $ and $ C\subseteq X,\ D\subseteq Y $, we write $ f''C $ for the image of $ C $ under $ f $, and $ f^{-1}[D] $ for the preimage of $ D $.
	
	Recall that a set $X$ is \emph{Dedekind-infinite} if there is an injection from $\omega$ into $X$, and is \emph{Dedekind-finite} otherwise.
	
	\begin{defn}
		A set $X \subseteq [\omega]^\omega$ is said to have the \emph{Ramsey property} if there exists an infinite set $H \in [\omega]^\omega$ such that either $[H]^\omega \subseteq X$ or $[H]^\omega \cap X = \emptyset$. We say that $H$ is homogeneous for $X$.
	\end{defn}
	
	\begin{defn}
		An \emph{almost disjoint family} is a collection $\mathcal A\subseteq[\omega]^\omega$ such that $A\cap B$ is finite for all distinct $A,B\in\mathcal A$. An almost disjoint family $\mathcal A$ is called \emph{maximal} $($or a \emph{MAD family}$)$ if for every $x\in[\omega]^\omega$, there exists $A\in\mathcal A$ such that $x\cap A$ is infinite.
	\end{defn}
	
	We note that a good pointclass $\Gamma$ is closed under intersections with analytic sets and has $A^{<\omega}\in\Gamma$ whenever $A\in\Gamma$.
	
	Now we can prove the main theorem.

	\begin{proof}[Proof of Theorem~\ref{thm:mad-ramsey}]
		If not, let $\mathcal A\in\Gamma$ be an infinite MAD family.
		
		\medskip
		\noindent\textbf{Case 1: $\mathcal A$ is Dedekind-infinite.}
		
		Fix a sequence $\{B_n:n\in\omega\}$ of distinct members of $\mathcal A$. For each $n\in\omega$, let $ A_n=B_n\setminus\bigcup_{i<n}B_i $. Then the $ A_n $ are infinite and pairwise disjoint.
		
		For $z\in\infsub{\omega}$, define
		\[
		\widetilde{z}
		=\{A_{z(2n)}(z(2n+1)):n\in\omega\}.
		\]
		Note that for all $ i\in\omega $ we have $ \widetilde{z}\cap B_i $ is finite.
		
		\begin{claim}\label{claim:dim1hide}
			For every $w\in\infsub{\omega}$, there exist $y\in\infsub{w}$ and $A\in\mathcal A$ such that $\widetilde{y}\subseteq A$.
		\end{claim}
		
		\begin{claimproof}{\ref{claim:dim1hide}}
			Since $\mathcal A$ is a MAD family and $\widetilde{w}$ is infinite, there exists $A\in\mathcal A$ such that $\widetilde{w}\cap A$ is infinite. Increasingly enumerate $ \widetilde{w}\cap A  $ as
			\[
			\{A_{w(2n_k)}(w(2n_k+1)): k\in\omega\}.
			\]
			Define
			\[
			y=\{w(2n_k):k\in\omega\}
			\cup\{w(2n_k+1):k\in\omega\}.
			\]
			Then
			\[
			\widetilde{y}
			=\{A_{w(2n_k)}(w(2n_k+1)):k\in\omega\}
			=\widetilde{w}\cap A\subseteq A,
			\]
			as required.
		\end{claimproof}
		
		\begin{claim}\label{claim:dim1split}
			For every $w\in\infsub{\omega}$, there exist $y\in\infsub{w}$ and $A\in\mathcal A$ such that both $A\cap\widetilde{y}$ and $\widetilde{y}\setminus A$ are infinite.
		\end{claim}
		
		\begin{claimproof}{\ref{claim:dim1split}}
			By the previous claim, we can assume that $\widetilde{w}\subseteq A$ for some $A\in\mathcal A$. Note that we have $A\neq B_i$ for every $i\in\omega$, so $A\cap A_i\subseteq A\cap B_i$ is finite for every $i\in\omega$.
			
			Let $w_e=\{w(2n):n\in\omega\}$. We construct $y\subseteq w$ by giving a recursive definition of its increasing enumeration. For convenience, let $y(-1)=-1$. Assuming that we have already defined $y(j)$ for $j<i$, where $i$ is divisible by $4$, we define the next four elements as follows:
			\begin{itemize}
				\item Define
				$y(i)=\min\{u\in w_e:u>y(i-1)\}$.
				
				\item Define
				\[
				y(i+1)=\min\{u\in w:u>y(i)
				\text{ and }A_{y(i)}(u)\notin A\}.
				\]
				This is possible since $A_{y(i)}\cap A$ is finite.
				
				\item Define
				$y(i+2)=\min\{u\in w_e:u>y(i+1)\}$.
				
				\item Define $y(i+3)=\min\{u\in w:u>y(i+2)\}$. Since $y(i+2)\in w_e$ and $\widetilde{w}\subseteq A$, we have $A_{y(i+2)}(y(i+3))\in A$.
			\end{itemize}
			Thus both $A\cap\widetilde{y}$ and
			$\widetilde{y}\setminus A$ are infinite.
		\end{claimproof}
		
		Let us define
		\[
		P=\{z\in\infsub{\omega}:
		\exists A\in\mathcal A\
		(\widetilde{z}\subseteq A)\}.
		\]
		Since the map $z\mapsto\widetilde{z}$ is continuous, the relation $\widetilde{z}\subseteq A$ is closed. By the closure properties of $\Gamma$, we have $P\in\Gamma$.
		
		We now show that $P$ is not Ramsey.
		\begin{itemize}
			\item Suppose that there exists $w\in\infsub{\omega}$ such that $\infsub{w}\cap P=\emptyset$. By Claim~\ref{claim:dim1hide}, there are $y\in\infsub{w}$ and $A\in\mathcal A$ with $\widetilde{y}\subseteq A$, i.e., $y\in P$, a contradiction.
			\item Suppose that there exists $w\in\infsub{\omega}$ such that $\infsub{w}\subseteq P$. By Claim~\ref{claim:dim1split}, there exist $y\in\infsub{w}$ and $A\in\mathcal A$ such that both $A\cap\widetilde{y}$ and $\widetilde{y}\setminus A$ are infinite. Since $y\in P$, there exists $B\in\mathcal A$ such that $\widetilde{y}\subseteq B$. So $A\cap B$ is infinite, which implies that $A=B$. This means that $\widetilde{y}\subseteq A$, a contradiction.
		\end{itemize}
		Therefore $P$ is a non-Ramsey set in $\Gamma$, again a contradiction.

		\medskip
	\noindent\textbf{Case 2: $\mathcal A$ is Dedekind-finite.}
	
	Note that for any finite $F\subseteq\mathcal A$, the set $\omega\setminus\bigcup F$ is infinite. For any finite subset $F\subseteq[\omega]^\omega$ with $\omega\setminus\bigcup F$ infinite, let $e_F$ be the increasing enumeration of $\omega\setminus\bigcup F$.
	
	For $n\in\omega$, define 
	\[
	D_n=\left\{x\in\infsub{\omega}:\exists \left(A_0,\ldots,A_{n-1}\right)\in{\mathcal A}^n \ 
	\forall j<n \ \ 
	e_{\{A_i:i<j\}} {}''x\subseteq^* A_j \right\}.
	\]
	
	By the definition of $ D_n $, any tuple $ \left(A_0,\ldots,A_{n-1}\right) $ witnessing $ x\in D_n $ has pairwise distinct entries; for $ n=0 $ the empty tuple is the unique witness. Also, for every $n\in\omega$, we have $D_{n+1}\subseteq D_n$. Moreover, if $x\in D_n$ and $y\in\infsub{\omega}$ has $y\subseteq^*x$, then $y\in D_n$.
	
	\begin{claim}\label{claim:df-unique}
		For every $n\in\omega$ and $x\in D_n$, the tuple $(A_0,\ldots,A_{n-1})$ witnessing $x\in D_n$ is unique.
	\end{claim}
	
	\begin{claimproof}{\ref{claim:df-unique}}
		Suppose $(A_i)_{i<n}$ and $(B_i)_{i<n}$ are two distinct tuples witnessing $x\in D_n$. Let $j<n$ be least such that $A_j\neq B_j$, and let $ F=\{A_i:i<j\}=\{B_i:i<j\} $. Then we have $ {e_F}''x\subseteq^*(A_j\cap B_j) $, contradicting almost disjointness.
	\end{claimproof}
	
	\begin{claim}\label{claim:df-intersection}
		$\displaystyle\bigcap_{n\in\omega}D_n=\emptyset$.
	\end{claim}
	
	\begin{claimproof}{\ref{claim:df-intersection}}
		Suppose $x\in\bigcap_{n\in\omega}D_n$. For any $ m\in\omega $ there is a unique $ s_m \in \idl{A}^{m+1} $ witnessing $ x\in D_{m+1} $. For any $ 0\leq m <n$, as $ s_n \upharpoonright (m+1) $ also witnesses $ x\in D_{m+1} $, we have $ s_n \upharpoonright (m+1)  = s_m$. Define a function $ f:\omega \to \idl{A} $ by $ f(n) = s_n(n) $. Since the entries of $s_n$ are pairwise distinct, the function $ f $ is injective. This contradicts the Dedekind-finiteness.
	\end{claimproof}
	
	\begin{claim}\label{claim:df-density}
		For every $n\in\omega$ and $H\in\infsub{\omega}$, we have $\infsub{H}\cap D_n\neq\emptyset$.
	\end{claim}
	
	\begin{claimproof}{\ref{claim:df-density}}
		We argue by induction on $n$.
		\begin{itemize}
			\item For $n=0$, we have $D_0=\infsub{\omega}$.
			\item Suppose the claim holds for $n$, and fix $H\in\infsub{\omega}$. To prove the claim for $ n+1 $, fix $x\in\infsub{H}\cap D_n$ with the witness $(A_0,\ldots,A_{n-1})$, and let $F=\{A_i:i<n\}$. By maximality, there is $A\in\mathcal A\setminus F$ such that $e_F{}''x\cap A$ is infinite. Thus $y=x\cap e_F^{-1}[A]$ is also infinite. We have that $(A_0,\ldots,A_{n-1}, A)$ witnesses $ y\in D_{n+1} $. So $y\in\infsub{H}\cap D_{n+1}$.
		\end{itemize}
	\end{claimproof}
	
	Define
	\[
	E=\{x\in\infsub{\omega}:x\in D_{\min x+1}\}.
	\]
	
	Equivalently,
	\[
	E=\left\{x\in[\omega]^\omega:
	\exists s\in\mathcal A^{<\omega}\ 
	\left[
	|s|=\min x+1\ \wedge\
	\forall j<|s|\
	e_{\{s(i):i<j\}}{}''x\subseteq^*s(j)
	\right]\right\}.
	\]
	Since the relation inside the brackets is analytic and $\mathcal A^{<\omega}\in\Gamma$, the closure properties of $\Gamma$ imply that $E\in\Gamma$.
	
	We show that $E$ is not Ramsey.
	
	\begin{itemize}
		\item Suppose there is $H\in\infsub{\omega}$ with $\infsub{H}\cap E=\emptyset$. Let $ k=\min H $. By Claim~\ref{claim:df-density}, take $x\in\infsub{H\setminus\{k\}}\cap D_{k+1}$. As $\{k\}\cup x\in D_{k+1}$, we have $\{k\}\cup x\in E$, a contradiction.
		\item Suppose there is $H\in\infsub{\omega}$ such that $\infsub{H}\subseteq E$. For every $ n\in\omega $, we have $ H\setminus n \in D_{\min(H\setminus n)+1} \subseteq D_{n+1} $, which implies $ H\in D_{n+1} $. Thus $H\in\bigcap_{n\in\omega}D_{n}$, contradicting Claim~\ref{claim:df-intersection}.
	\end{itemize}
	
	Therefore $E$ is a non-Ramsey set in $\Gamma$, a contradiction.
		
	This completes the proof.
	\end{proof}
	
	Theorem \ref{thm:mad-ramsey} has several classical consequences. In particular, it gives an alternative proof that there are no analytic infinite MAD families.
	
\begin{cor}[\cite{mathias1977happy}]
	There are no analytic infinite MAD families.
\end{cor}

\begin{proof}
	By absoluteness, in $ \zfaxiom $ every analytic set is Ramsey. Apply Theorem~\ref{thm:mad-ramsey} with $\Gamma=\boldsymbol{\Sigma}^1_1$.
\end{proof}
	
	As a further application, we recover Törnquist's theorem on Solovay's model.
	
	\begin{cor}[\cite{tornquist2018definability}]
		There are no infinite MAD families in Solovay's model.
	\end{cor}
	\begin{proof}
		In Solovay's model \cite{solovay1970model}, every set of reals in this model has the Ramsey property \cite{mathias1969generalization}. Apply Theorem \ref{thm:mad-ramsey} where $\Gamma$ is the pointclass of all sets.
	\end{proof}
	
	Assuming the \emph{Projective Determinacy} $ \mathsf{PD} $, we also obtain:
	
	\begin{cor}[\cite{neeman2018happy}]
		$\mathsf{PD}$ implies there are no projective infinite MAD families.
	\end{cor}
	\begin{proof}
		Under $\mathsf{PD}$, every projective set is Ramsey \cite{harrington1981determinacy}. Suppose $ \idl{A} $ is a projective infinite MAD family. Find $ n\in\omega $ such that $ \idl{A}\in{\boldsymbol \Sigma}^1_n $. Apply Theorem \ref{thm:mad-ramsey} where $\Gamma = {\boldsymbol \Sigma}^1_n$.
	\end{proof}
	
	\section*{Acknowledgments}
	
	The authors are grateful to Renling Jin, David Schrittesser, and Hang Zhang for many valuable comments and discussions that helped improve this paper. Jialiang He was supported by the National Natural Science Foundation of China (Grant No.~12671002).
	
	\bibliographystyle{plain} \bibliography{mybibtex}
	
\end{document}